\documentclass[12pt]{article}

\usepackage{amssymb} \usepackage{amsfonts}

\begin{document}








\newcommand{\calR} {{\bf {\cal R}}}
\newcommand{\CalR} {{\bf {\cal R}}^\ast }

\newcommand{\R} {{\mathbb R}}
\newcommand{\Q} {{\mathbb Q}}
\newcommand{\C} {{\mathbb C}}
\newcommand{\F} {{\mathbb F}}
\newcommand{\Z} {{\mathbb Z}}
\newcommand{\N} {{\mathbb N}}
\newcommand{\card} {{\# } }
\newcommand{\mod}{\; {\rm mod}\; }
\newcommand{\cliff}{{\rm Cliff}}
\newcommand{\cliffc}{{\rm Cliff}_{\C }}
\newcommand{\fcm}{{\cal F}_{CHD}}
\let\lra=\Leftrightarrow
\let\nd=\noindent

\newcommand{\sig}{{\Sigma }_{2k}}
\newcommand{\sgn}{{\rm sgn\hspace{0.05cm}}}
\newcommand{\calB}{{\cal B}}

\def\fracsmall#1#2{{\textstyle {\frac {#1}{#2} }}}

\let\la=\langle
\let\ra=\rangle
\def\Qed#1{\mathop{\mkern0.5\thinmuskip\vbox
{\hrule\hbox{\vrule\hskip#1\vrule height#1
width 0pt\vrule}\hrule}\mkern0.5\thinmuskip}}
\def\qed{$\Qed{8pt}$}


\newtheorem{theo}{Theorem}
\newtheorem{theorem}{Theorem}
\newtheorem{prop}{Proposition}
\newtheorem{lemma}{Lemma}
\newtheorem{corollary}{Corollary}
\newtheorem{examp}{Example}
\newtheorem{remark}{Remark}


\title{Flat manifolds, harmonic spinors, and eta invariants}
\author{{M. Sadowski, A. Szczepa\'nski}
\thanks{Supported by University of Gda\'nsk grant number BW - 5100-5-0319-3}}
\date{\today}
\maketitle

\begin{abstract} 
The aim of this paper is to calculate the eta invariants and the 
dimensions of the spaces of harmonic spinors of an 
infinite family of closed flat manifolds 
$\fcm.$ It consists of some flat manifolds $M$
with cyclic holonomy groups. 
If  $M \in \fcm ,$  then we give explicit formulas for 
$\eta (M) $ and  ${\mathfrak h}(M). $ 
The are expressed in terms of solutions of 
appropriate congruences in $\{ -1,1 \}^{ \left [\frac {n-1}{2} \right ]}.$  
As an application we investigate the integrability of some $\eta$ invariants of
$\fcm$-manifolds.
\end{abstract}

\vspace{0.2cm} 
\section{Introduction}
In this paper we consider Dirac operators on
an infinite family $\fcm $  of closed flat manifolds.
It consists of flat manifolds $M$
with cyclic holonomy groups of odd order equal to
the  dimension of $M.$
The family $\fcm $ is particularly simple
and the investigation of different properties
of multidimensional flat manifolds should start with
the investigation of them in this particular case.
Some $\fcm$ manifolds arises in the classification
of flat manifolds whose holonomy groups have prime order (cf. \cite{Ch}).
We describe the eta invariants of the Dirac operators
arising from different spin structures and we
give necessary and sufficient conditions
of the existence of nontrivial harmonic spinors.
The methods used here extends that used in \cite{Pf}.
We apply them to much wider class of manifolds
and we consider related general questions.

\vspace{0.2cm}
To formulate the main results we need some definitions.
Let $n=2k+1$ be an odd number, and let
$a_1,...,a_n$ be a  basis of $\R ^n$.
Consider the linear map
 $A :\R ^n \to \R ^n$ such that
$A(a_j)=a_{j+1}$ for $j<n-1$,
$A(a_{n-1})= -a_1 -...-a_{n-1}, $ and
$A(a_n)=a_n.$
Let $a=\frac {1}{n} a_n$ and let
$g(x)=A(x) +a .$
An $n$-dimensional flat manifold $M \in \fcm  $
can be written as $\R ^n /\Gamma ,$
where $\Gamma = \la a_1,...,a_{n-1}, g \ra .$
The linear part $A$ of $g$ has two lifts
$\alpha _{+}, \alpha _{-} \in Spin (n) $ such that
$ \alpha _{+}^n =id $ and $\alpha _{-}^n=-id $ (see Section
\ref{sec-spin-struct-fcm}).
This defines two spin structures on $M.$

To formulate the result describing
$\eta _{M^{n}}(0)$
for $M^{n}\in  \fcm  $  we need some
combinatorial invariants.
It is known that $\eta _{M^n}(0)=0 $ if $k$ is even
(cf. \cite[p. 61]{aps}) so we consider the case when $k$ is odd.
Let
$$ c(k)=\cases { 0 &  if $\frac {k(k+1)}{2}$ is even \cr
\frac {1}{2} &  if $\frac {k(k+1)}{2}$ is odd \cr
}
.$$
For every $\epsilon =(\epsilon _1,...,\epsilon _k) \in \{ -1,1 \} ^k $
consider
$$\mu _\epsilon = \sum _{j=1}^k \epsilon _j j
\quad {\rm and} \quad
\nu (\epsilon ) =\epsilon _1 \cdots \epsilon _k .
$$
Let ${\cal D}_{+} = \{ \epsilon \in \{-1,1 \} ^k :  \nu (\epsilon ) =
\epsilon_{1}\epsilon_{2}...\epsilon_{k} = 1 \}, $
let $r\in \{ 0,...,n-1 \},$  let
$$  A_{r}^{+} = 2 \card \{ \epsilon \in {\cal D}_{+}   :\,
\frac {\mu _\epsilon }{2} +c(k)n \equiv  r \mod (n)  \} $$
in the case of $\alpha  _{+} ,$
and let
$$  A_{r}^{-} = 2 \card \{ \epsilon \in {\cal D}_{+}  :\,
\frac {\mu _\epsilon }{2} +c(k)n +k \equiv  r \mod (n) \} $$
in the case of $\alpha  _{-} .$
The numbers $A_{r}^{\pm}$ are well defined
(cf. Remark \ref{rem-congruence-well-defined}).

\begin{theo}\label{thm-calcul-eta}\hspace{-0.2cm}{\bf .} \hspace{0.05cm}
Let $k$ be an odd positive integer and let $n=2k+1.$
If $M^n \in \fcm ,$ then
$$ \eta _{M^{n},\alpha _{+}} (0) =
\sum _{r =1 }^{n-1} A_{r}^{+} \left (1- \frac {2r}{n } \right )
, \leqno (1) $$
$$ \eta _{M^{n},\alpha _{-}}(0) =
\sum _{r =0 }^{n-1} A_{r}^{-} \left (1- \frac {2r+1}{n } \right )
 .  \leqno (2) $$

\end{theo}

Applying Theorem 1 we prove that some $\eta$-invarints of $\fcm$-manifolds
are integral (Corollary 1) and that $\eta_{M,\alpha_{+}}-\eta_{M,\alpha_{-}} \in 2\Z$
(Corollary 2).  
\noindent
Let ${\mathfrak h}(V)$ be the dimension of the vector space of harmonic
spinors.
\vskip 1mm
\noindent
Now we state another result of the paper.

\begin{prop}\label{thm-harmonic-spinors}\hspace{-0.2cm}{\bf .}
\hspace{0.05cm}
Let $k$ be a positive integer and let $n=2k+1.$
If  $M\in \fcm ,$ then

\vspace{0.1cm}
\noindent
{\bf a)}
${\mathfrak h}(M,\alpha _{+}) >0 $  if and only if
$n\geq 5.$

\vspace{0.1cm}
\noindent
{\bf b)}
${\mathfrak h}(M,\alpha _{-}) =0. $
\end{prop}

\vspace{0.2cm}
The spectra of the Dirac operators on flat tori were described
in \cite{TF-tori}.
The spectra of the Dirac operators
on closed $3$-dimensional
flat manifolds and their eta invariants
were calculated in \cite{Pf}.
We should also mention about (\cite{MP}) where the authors consider spin structure
and the Dirac operators on flat manifolds with $\Z_{p},$ (p-prime number), and non-cyclic holonomy.

\vspace{0.2cm}
Throughout this paper the following notation will be used.
If $G$ is a group and $g_1,...,g_l \in G$, then
$\la g_1,...,g_l \ra $ is the subgroup of $G$ generated by
$g_1,...,g_l.$
The symbol $X^G$ stands for the set of the fixed points of an
action of $G$ on $X$.
For every $g\in G$, $X^g =\{ x\in X: gx=x \} .$
By $\Gamma $ (or $\Gamma _n$) we denote the deck group of a
closed flat manifold $M,$
by $h$ the holonomy homomorphism of $M$, and by $\widehat h$
its lift to $Spin(n).$
The standard epimorphism from $Spin(n)$ to $SO(n)$
will be denoted by $\lambda $ (cf. Section \ref{sec-spin-struct-fcm}).
The letter $\Gamma _0$ stands for the maximal abelian subgroup
of $\Gamma $ consisting of all translations belonging to
$\Gamma ,$ (cf. \cite{Ch} and \cite{wolf}).
By $a_1,...,a_n$ we usually denote a basis of $\Gamma _0.$
The subspace of a vector space spanned by vectors
$v_1,...,v_l$ will be denoted by ${\rm Span} [v_1,...,v_l].$
The symbols $\alpha _{+},$
$\alpha _{-},$  ${\mathfrak h}(M)$,
$c(k),$ $\mu _\epsilon $, $\nu (\epsilon ) $,  and $A_r$ were defined above.
The cyclic group $\la A \ra $ will be denoted by $G.$
\vskip 1mm
\noindent
We would like to thank Andrzej Weber for helpful conversations.
We are grateful to Roberto
Miatello for correcting a mistake in an earlier version of the paper and to 
Bernd  Ammann for pointing out a typographic error.


\section{Spin structures on $\fcm $-manifolds\\ and Dirac
operators}\label{sec-spin-struct-fcm}

Let $k \in \N \cup \{ 0\} $ and let $n=2k+1.$
Let $\Gamma$ be as in the introduction,
and let
$\la , \ra  ^\ast  $ be an $A$-invariant
scalar product in $\R ^n.$
From definition (cf. \cite{Ch} and \cite{wolf})
$ M=\R ^n /\Gamma $ is a closed, orientable, flat manifold.
Moreover the eigenvalues of the generator $A$ of the holonomy group of $M$ are equal to
$e^{\frac {2\pi i j}{n} } , j=1,...,n.$ In fact,
for every $j=2,...,n-1,$ consider the
$(j\times j)$-matrix
$$
M_j=
\left [
\begin{array}{lllll}
0 & 0 & \dots & 0 & -1 \\
1 & 0 & \dots & 0 & -1 \\
0 & 1 & \dots & 0 & -1 \\
\vdots & \vdots &  & \vdots & \vdots \\
0 & 0 & \dots & 1 & -1 \\
\end{array}
\right ]
.$$
Let $M_j(z)=M_j- zI$. Then
$$ \det (A-zI)= (1-z)\det M_{n-1}(z) .$$
Applying the Laplace expansion with respect to the first row
we have
$$ \det M_j(z) =-z \det M_{j-1}(z) +(-1) ^j .$$
Using this it is easy to check that
$ \det M_j(z) = (-1)^j \sum _{l=0}^j z^l .$
Hence
$$ \det (A-zI)= (1-z)\det M_{n-1}(z) =-z^n +1.$$
\vspace{0.2cm}
Let $e_1,...,e_n$ be an orthonormal basis in $(\R ^n, \la ,\ra ^\ast). $
Throughout the rest of the paper we shall always assume 
(cf. \cite[page 61]{aps} and \cite[Proposition 1.3]{Hitchin}) that:
\vskip 3mm
\noindent
(i)
$e_1,...,e_{n-1} \in {\rm Span}\, [a_1,...,a_{n-1}]$ and
$e_n = a_n,$

\vskip 1mm
\noindent
(ii) for every for $j\leq n-1  : \, $
$A(e_{2j-1} ) = \cos (2\pi j/n) e_{2j-1} +\sin (2\pi j/n)e_{2j} ,$
and
$A(e_{2j} ) = -\sin (2\pi j/n) e_{2j-1} +\cos (2\pi j/n)e_{2j} . $
\vskip 3mm
\noindent
Let $\cliff (n) $ be the Clifford algebra in $\R ^n$ and let
$\cliffc (n)$ be its complexification.
The group $Spin(n)$ is the set of products
$x_1\cdots x_{2r},$ where $r\in \N,$
and where $x_1,...,x_{2r}$ are the elements of the unit sphere in
$\R ^n.$
The standard covering map
$\lambda : Spin (n) \to SO(n) $ carries $y\in Spin (n)$
onto $\R ^n \ni x \to yxy^\ast, $
where $(e_{j_1}\cdots e_{j_s})^\ast = e_{j_s}\cdots e_{j_1}.$
A spin structure on an orientable flat manifold $M = \R ^n /\Gamma$
is determined by the lift $\widehat h :\Gamma  \to Spin (n)$ of
the holonomy homomorphism  $ h :\Gamma  \to SO (n)$.
Recall that  $h$ carries
$\gamma  \in \Gamma$ onto its linear part $h(\gamma ),$ (cf. \cite[Chapter III]{wolf}).
For $M\in \fcm $ we have
 $h(\Gamma ) =\la A \ra \cong \Z _n$
and  any lift $\widehat A$
of $A$ to $Spin (n)$ defines the lift
$\widehat h$ of $h$, given by the formulas
$ \widehat h (a_j) = 1$ for $j\leq n-1,$ $\widehat h(g)=\widehat A.$
In order to construct $\widehat A$
consider $\beta = \frac {\pi }{n},$
 $$ r_{j} = \cos(j\beta ) + e_{2j-1}e_{2j}\sin(j\beta ), $$
and  $\alpha = \prod _{j=1}^k r_j .$
Clearly $r_ir_j=r_jr_i $ for $i,j\in \{ 1,...,k \} .$
A direct calculation yields
$$
\lambda (r_j)(e_l)=
\cases {
 \cos (2j \beta ) e_{2j-1} + \sin (2j \beta )  e_{2j} &
for  $\; l=2j-1$ \cr
- \sin (2j \beta ) e_{2j-1} + \cos (2j \beta )  e_{2j} &
for  $\; l=2j$ \cr
 e_{l} &
for  $\; l \notin \{ 2j-1,2j \} $ \cr
} .
 $$
Using this it is easy to check that
 $$\alpha ^n = (-1) ^{\frac {k(k+1)}{2} } $$
and   $\lambda (\alpha )= A.$
Now we can define
$$ \alpha _{+} = (-1) ^{\frac {k(k+1)}{2} } \alpha ,
\quad
\alpha _{-} = -(-1) ^{\frac {k(k+1)}{2} } \alpha .$$
Since $n$ is odd,
$$ \alpha _{+}^n=1 \quad {\rm and} \quad
\alpha _{-}^n=-1. $$ 
We have.
\begin{lemma}\label{lem-first-homol-M}\hspace{-0.2cm}{\bf .} \hspace{0.05cm}
 $H_1(M,\Z ) \cong \Z \oplus H,$
where $H$ is a finite abelian group of odd order and
$H^1(M, \Z _2) \cong \Z _2.$
\end{lemma}

\vspace{0.1cm}
\noindent
{\bf Proof:}
The group $\Gamma_{0} =\la a_1,...,a_n \ra $ is the maximal
abelian subgroup of $\Gamma $ and the following sequence
$$ 0 \to \Gamma _0 \to \Gamma \to \la A \ra \to 1 $$
is exact (cf. \cite[Proposition 4.1]{Ch}, \cite[Theorem
3.2.9]{wolf}).
From \cite[Corollary 1.3]{HS} we have
$\dim_{\Q}(\Q \otimes H_1(\Gamma ,\Z ) = \dim_{\Q}(\Q \otimes \Gamma _0 ^A) = 1. $
Hence
 $H_1(M,\Z ) \cong \Z \oplus H$, where $H$ is a finite group.
According to \cite[Chapter 3]{brown}, there are homomorphisms
$res : H_\ast (M,\Z )\cong H_\ast (\Gamma ,\Z )  \to
H_\ast (\Gamma _0, \Z ) $ and
$cor : H_\ast ( \Gamma _0, \Z ) \to
H_\ast (M,\Z ) $ such that
$cor \circ res $ is the multiplication by $n.$
Since the group $H_\ast (\Gamma _0,Z )\cong H_\ast (T^n,\Z )$
is torsion free we have  $nH=0$. In particular, the order of $H$ is odd.
For the proof of the last statment we have
$H^1(M,\Z _2) \cong {\rm Hom}\, (H_1(M,\Z ), \Z _2 ) \cong
{\rm Hom}\, (\Z , \Z _2 )\cong \Z _2 .$
 \qed

\vspace{0.3cm}
Since $\alpha _{+} ,$  $\alpha _{-} $ are different lifts
of the holonomy homomorphism $h$ to $Spin(n)$,
the spin structures determined by them are different.
It is known that spin structures on $M$ correspond
to the elements of $H^1(M,\Z _2)$ (\cite[p.~40]{TF}).

\vskip 5mm
By \cite[Section 1.3]{TF}, the irreducible complex
$\cliffc (n)$-module $\Sigma _{2k}$ can be described as follows.
Consider
$$ g_1 =\left [
\matrix {
i & 0 \cr
0 & -i }
\right ],
\quad
g_2 =\left [
\matrix {
0 & i \cr
i & 0 }
\right ],
\quad
T =\left [
\matrix {
0 & -i \cr
i & 0 }
\right ]. $$
Let
$\sig = \underbrace {\C ^2 \otimes ... \otimes \C ^2 }_{k \; {\rm times}}$
and let
$\alpha (j) = \cases { 1 & if $j$ is odd \cr
2 & if $j$ is even \cr
} .$
Take an element $u=u_1 \otimes ...\otimes u_k $ of  $ \sig $
and the orthonormal basis $e_1,..,e_n$  considered above.
Then
$$ e_j u =(I\otimes ... \otimes I \otimes g_{\alpha (j)} \otimes
\underbrace { T \otimes ... \otimes T }
_{ [ \frac {j-1 }{2} ]  \; {\rm times} } )(u), $$
for $j \leq n-1$, and
$$ e_n u = i(T\otimes ... \otimes T) u .$$
A spin structure on $M$
determines
a complex spinor bundle $P\Sigma _{2k} $ with fiber
$\Sigma _{2k} $.
This bundle is the orbit space of
${\bf R}^n \times \Sigma _{2k} $ by the action of
$\Gamma $ given by
\begin{equation}\label{fibre}
\gamma (x,v)= (\gamma x, \widehat h(\gamma )v ) ,
\end{equation}
where $\gamma \in \Gamma, x \in {\bf R}^n$ and $v \in \Sigma_{2k}.$
Clearly
$$ \widehat h ( a_j) =1  \; {\rm for} \;  j\leq n-1 $$
and
$$ \widehat h ( g) = \alpha_{\pm} .$$
Since ${\rm Span}[e_1,...,e_{n-1}] =
{\rm Span}[a_1,...,a_{n-1}]$ and $a_n=e_n$ we conclude that
$$ \widehat h ( e_j) =1  \; {\rm for} \;  j\leq n-1 .$$
\noindent
Consider the covering $T^n=\R ^n/\Gamma _0 \to M.$
We have
$  \widehat h ( a_n) =
 \pm 1 .$
The lift  $P_T\Sigma _{2k}$  of $P\Sigma _{2k}$ to $T^n$
is the orbit space
$ (\R ^n \times \Sigma _{2k} )/\Gamma _0,$
where the action of $\Gamma _0$ on $ \R ^n \times \Sigma _{2k} $
is given by the formula (\ref{fibre}).

\vspace{0.2cm} To deal with the spectrum of the Dirac operator $D$ it
is convenient to describe it in terms of the spectrum
of $D^2$. We state without proofs some related results
of \cite{Pf} that will
be used later. Identify the parallel section $\R ^n \ni x \to
(x,v) \in \R ^n \times \Sigma _{2k}  $ with $v.$ Every section (spinor)
of the trivial bundle $\R  ^n \times \Sigma _{2k} $ (covering
our bundle $P\Sigma _{2k} $) can be written as a linear combination of $fv,$ where
$f\in C^{\infty }( \R ^n, \C) $ and $v$ is a parallel section.
 Take the coordinate system $x_1,...,x_n$ determined by
$e_1,...,e_n$.
Since $v$ is parallel,
\begin{equation}\label{dirac}
D(fv) =\sum _j e_j \nabla _{e_j} (fv)=
\sum _j e_j \left ( \fracsmall
{\partial }{\partial x_j} (f)v + f\nabla _{e_j} v
\right ) =
\sum _j e_j \fracsmall {\partial f}{\partial x_j} v .
\end{equation}
Let $  \Gamma _0 ^\ast$ be the dual lattice of $\Gamma_{0}.$
Let $\calB $ be $\Gamma _0 ^\ast $ in the case of ${\alpha }_{+}$
and  $\Gamma_{0} ^\ast +\frac {1}{2} e_n$ in the case of  ${\alpha }_{-}$.
The action of $g$ on on the set of sections
of $\R ^n \times \sig $, induced by the action of $g$ on $\R ^n,$
is given by the formula
\begin{equation}\label{spinor}
g(\phi)(x)= \widehat h (g) \phi(g^{-1}x) ,
\end{equation}
where $\phi$ is a spinor on $\R^{n}.$

Consider $f_b(x)=e^{2\pi i \langle b, x \rangle }.$
By immediate calculation or following (\cite{Pf}) we have
\begin{equation}\label{dsq}
D^2(f_bv) =4\pi ^2 \vert \vert b \vert \vert ^2  f_b v .
\end{equation}
Hence the sections $f_bv ,b \in \calB   , $
$v\in \sig $,   correspond to eigenvectors of
 $D$ on $T^n ,$
and
the elements of $\{ f_bv :v \in \sig ,b\in \calB \}^g $
 correspond to
eigenvectors of  $D$ on $M .$

For $b\in \calB,$ let us denote the corresponding $D^{2}$-eigenspace
by $E_{b}(D^{2}) = {\rm Span}\{f_{b}v : v\in \sig \}.$
We have the decomposition
$E_b(D^2)=E_{b+}(D)\oplus E_{b-}(D)$,
where
$$ E_{b\pm }(D) = \{ p\in E_b(D^2) :\,
Dp=\pm 2\pi \vert \vert b \vert \vert \, p \}. $$
Since
\begin{equation}\label{action}
(f_b \circ g^{-1} )(x) =e^{-2\pi i \la A(b),a \ra } f_{A(b)}(x) 
\end{equation}
we have $ AE_b \subseteq E_{A(b) } ,$ (cf. \cite[Lemma 4.1]{Pf}).
Denote  $<A>$ by $G.$
\vskip 1mm
\noindent
Let $\calB _{Sym}=\{b\in \calB : \# G(b) = \#G   \} $, 
$\calB _{Pas}=\{b\in \calB  : \# G(b) <\#G \}$ and
\begin{equation}\label{dec}
D_S =D\vert _
{[\bigoplus _{b \in {\calB _{Sym} }} E_b(D^2)]^g } , \quad 
D_{Pas} =D\vert _
{[\bigoplus _{b \in {\calB _{Pas} }}E_b(D^2) ]^g } .
\end{equation}
Clearly $\calB $ is the disjoint union of
$\calB _{Sym}$ and $\calB _{Pas} $ and
the Dirac operator $D$ on $M$ can be identified with
$ D_{S} \oplus D_{Pas} .$
If $b\in \calB_{Sym}$ and
$$V_b^{\pm}=\bigoplus \limits _{h\in G} E_{h(b\pm)}(D^{2})$$
then $\dim (V_b ^{\pm} )^g =\dim E_b ^{\pm }(D) =2^{k-1}$
(cf. \cite[Theorem 4.2, Corollary 4.3]{Pf}).

\section{Eta invariants of $\fcm $-manifolds}\label{sec-eta-for-fcm}

The aim of this section is to prove Theorem \ref{thm-calcul-eta}.
Recall that the $\eta $-invariant of the Dirac operator on a
closed spin manifold $M$ is defined as follows. As $D$ is elliptic
formally self adjoint, it has discrete real spectrum and the
series $ \sum _{\lambda \neq 0} \sgn \lambda \vert \lambda \vert
^{-z}$ converges for $z\in \C$ with ${\rm Re} (z)$ sufficiently
large (\cite[Theorem 3.10]{aps}). Here summation is taken over all
nonzero eigenvalues $\lambda $ of $D$, each eigenvalue being
repeated according to its multiplicity. The function $z\to  \sum
_{\lambda \neq 0} \sgn \lambda \vert \lambda \vert ^{-z}$ can be
extended to a meromorphic function $\eta _M$ in the whole complex
plane such that $0$ is a regular point of $\eta _M$ (\cite[Theorem
3.10]{aps}). The $eta $-{\it invariant} of $M$ is $\eta _M (0).$

\vspace{0.2cm}
Define an endomorphism $\rho _1$ of $\C ^2 $ by the formula
$$ \rho _1(u) = \cos \beta \, u +\sin \beta \, g_1g_2 \, u .$$
The matrix of $\rho _1$ is equal to
$$
\cos \beta \, I + \sin \beta
\left [ \begin{array}{rr}
0 & -1 \\
1 & 0 \\
\end{array}
\right ]
=
\left [ \begin{array}{rr}
\cos \beta & -\sin \beta \\
\sin \beta & \cos \beta \\
\end{array}
\right ]
$$
so that the matrix of $\rho _1^j $ is equal to
$$
\cos (\beta j) \, I + \sin (\beta j)
\left [ \begin{array}{rr}
0 & -1 \\
1 & 0 \\
\end{array}
\right ]
.
$$
The following lemma is crucial.

\begin{lemma}\label{lem-eigenvecs-of-alpha}\hspace{-0.2cm}{\bf .}
\hspace{0.05cm}
Let $w_{+1} =(1,-i),$ let $w_{-1} =(1,i),$
let $\epsilon =(\epsilon _1,...,\epsilon _k ) \in
\{ -1,1 \}^k, $
and let
$v_\epsilon =w_{\epsilon _1} \otimes ... \otimes w_{\epsilon _k}.$
Let $\beta =\frac {\pi }{n} $
and let $\mu _\epsilon $
be as in Theorem \ref{thm-calcul-eta}.
Take $u=u_1 \otimes .... \otimes u_k \in \sig .$
Then

\vspace{0.1cm}
\noindent
{\bf a)} $\alpha u = \rho _1u_1 \otimes .... \otimes \rho _1^k u_k, $

\vspace{0.1cm}
\noindent
{\bf b)} $\alpha e_n u = e_n \alpha u,$

\vspace{0.1cm}
\noindent
{\bf c)} $ \rho _1 (w_{\pm 1})=e^{\pm i \beta } w_{\pm 1} ,$

\vspace{0.1cm}
\noindent
{\bf d)} $\alpha v_\epsilon = e^{i \beta \mu _\epsilon}v_\epsilon $
and $\{ v_\epsilon : \epsilon \in \{-1,1 \}^k \} $
is a basis of $\sig ,$

\vspace{0.1cm}
\noindent
{\bf e)} $ e_n v_\epsilon  = - i \nu (\epsilon )  v_\epsilon .$
\end{lemma}

\vspace{0.1cm}
\noindent
{\bf Proof.}
{\bf a)}
Since $T^2=id,$
$$ e_{2j-1}e_{2j}(u_1\otimes ... \otimes u_j \otimes ... \otimes u_k ) =
u_1\otimes ... \otimes  g_1g_2 (u_j )\otimes ... \otimes u_k
 $$
and consequently
$$ r_j( u_1\otimes ... \otimes u_j \otimes ... \otimes u_k ) =
u_1\otimes ... \otimes \rho _1 ^j (u_j)  \otimes ... \otimes u_k
.$$
Hence
$$ \alpha (u_1\otimes ... \otimes u_k ) =
(r_1\cdots r_k)(u_1\otimes ... \otimes u_k )
= \rho _1 (u_1)\otimes ... \otimes \rho _1 ^k (u_k ) .$$

\vspace{0.2cm}
\noindent
{\bf b)}
For $j\leq k$ we have
$ e_{2j-1}e_{2j} e_n = e_n e_{2j-1}e_{2j} $ so that
$r _j e_n =e_n r_j .$

\vspace{0.2cm}
\noindent
{\bf c)} is obvious.

\vspace{0.2cm}
\noindent
{\bf d)}
By c),
$r_1(w_{\epsilon _j}) =e^{i\beta \epsilon _j } w_{\epsilon _j}.$
Hence
$$
\alpha (v_{\epsilon }) =
r_1(w_{\epsilon _1 }) \otimes ... \otimes
r_1^k(w_{\epsilon _k}) =
e^{i\beta \mu _\epsilon } v_\epsilon .$$
Since $\card \{ v_\epsilon :\, \epsilon \in \{ -1,1 \} ^k \} =2^k =
\dim \sig $ and the vectors $v_\epsilon $ are linearly independent,
they form a basis of $\sig .$

\vspace{0.2cm}
\noindent
{\bf e)}
We have
$T(w_1)=-w_1 $  and $T(w_{-1})=w_{-1}.$
It follows that
$$ e_n v_{\epsilon }= i T w_{\epsilon _1}
\otimes ... \otimes
T w_{\epsilon _k} =
i  (-1)^{ \card \{ j\in \{1,..,k \}:  \, \epsilon _j =1 \}  } v_{\epsilon }
= - i \nu (\epsilon )  v_\epsilon
 .$$
This finishes the proof of Lemma 2. \qed

\vskip 2mm
\noindent

Let ${\cal E}(\lambda , D_{Pas} )$ be the eigenspace of $\lambda $ for
 $D_{Pas}$ on $M$. 
 From the definition and (\ref{dsq}), (\ref{action}) it is easy to see that
 $\lambda = 2\pi l$ in the case of $\alpha_{+}$ and $\lambda = 2\pi (l+\frac{1}{2})$ in the case
of $\alpha_{-},$ where $l\in \Z.$
In the case $\alpha_{+}$, ${\cal B}_{Pas} = \{le_{n} : l\in  \Z\}$
and we have
\begin{equation}\label{basis}
D(f_{le_n} v_\epsilon )=
- \fracsmall {\partial  }{\partial x_n } (e^{2\pi i \la le_n , x \ra } )
i \nu  (\epsilon)  v_\epsilon =
 \nu (\epsilon ) 2\pi l f_{le_n} v_\epsilon .
 \end{equation}

\vskip0.2cm
\noindent
Hence 
${\cal E}(2\pi l , D_{Pas} ) = 
{\rm Span}[f_{\nu (\epsilon ) le_n }v_\epsilon : \epsilon \in \{-1,1 \}^{k}]^{g}.$
Similar formulas are also true for $\alpha_{-},$ where ${\cal B}_{Pas} = \{(l+\frac{1}{2})e_{n} : l\in \Z\}.$
\vskip0.2cm
Now we are able to describe the spectrum of $D_{Pas}.$

\begin{prop}\label{prop-spectrum-Dpas}\hspace{-0.2cm}{\bf .} \hspace{0.05cm}
Let $n$, $k,$ $M$, $\mu _\epsilon ,$
$\nu (\epsilon ) ,$ and $c(k)$ be as in Theorem \ref{thm-calcul-eta}.
Let $b^{+} = le_{n}$ and $b^{-} = (l+ \frac {1}{2})e_{n},$ where $l\in \Z.$
\vskip 1mm
\noindent
{\bf a)}
If the spin structure is given by
$\alpha _{+}  $,   then
 $$ {\cal E} (2\pi l ,D_{Pas} ) =
{\rm Span} [  f_{\nu (\epsilon ) b^{+}}v_\epsilon  :\;
\frac {\mu _\epsilon }{2} +c(k)n
\equiv  \nu (\epsilon ) l\; \mod (n) ].$$

\vskip 1mm
\noindent
{\bf b)}
If the spin structure is given by
$ \alpha _{-}$,
 then
$$ {\cal E} (2\pi(l+\frac {1}{2}) ,D_{Pas} ) =
{\rm Span} [
 f_{\nu (\epsilon ) b^{-}}v_\epsilon :\;
\frac {\mu _\epsilon }{2} +c(k)n +  \frac {n-1}{2}
\equiv  \nu (\epsilon ) l\; \mod (n) ]. $$
\end{prop}

\begin{remark}\label{rem-congruence-well-defined}\hspace{-0.2cm}{\bf .}
\hspace{0.05cm}
\begin{rm}
Since $\epsilon _j - 1 $ are even,
the difference
$\mu _\epsilon -k(k+1)/2 =
\sum _{j=1}^k \epsilon _j j -
\sum _{j=1}^k j $
is divisible by $2.$
Using this and the definition of $c(k)$
it is easy to see that $\mu _\epsilon /2 + c(k)n$
is an integer.
\end{rm}
\end{remark}

\vspace{0.2cm}
\noindent
{\bf Proof of Proposition \ref{prop-spectrum-Dpas}.} {\bf a)}
From the definitions of $\alpha _{+} $ and $c(k)$
it follows that
$\alpha _{+} = (-1) ^{2c(k)} \alpha .$
We have
$$g(f_{l e_n}(x)v_\epsilon ) =
f_{l e_n}(g^{-1}x) \alpha _{+} v_\epsilon  =
e^{- \frac {2\pi i l} {n} } f_{l e_n}(x) (-1)^{2c(k)} \alpha v_\epsilon $$
$$ = e^{-\frac {2\pi i} {n} (l - \fracsmall {\mu _\epsilon }{2} -c(k)n)  }
f_{l e_n}(x)  v_\epsilon .$$
By the above one gets the required conditions.

\vspace{0.1cm}
\noindent
{\bf b)}
In the case of $\alpha _{-}$ 
the eigenvectors of $D_{Pas }$ on $T^n$ can be written as
$f_{(l+\frac {1}{2})e_n} v$ for $l\in \Z, v\in \sig .$
We have
$$ g f_{(l+\frac {1}{2})e_n}(x) v_\epsilon =
e^{-\frac {2\pi i}{n} (l-c(k)n-\frac {n-1}{2} -\frac {\mu _\epsilon }{2}) }
f_{(l+\frac {1}{2})e_n}(x) v_\epsilon .$$
Hence $f_{({b^{-}})} v_\epsilon $ is $g$-equivariant
if and only if
$l \in n \Z + c(k)n+\frac {n-1}{2} + \frac {\mu _\epsilon }{2}.$
The rest of the argument is the same as in a).
This finishes the proof of
\hbox{Proposition \ref{prop-spectrum-Dpas}. \qed }

\begin{lemma}\label{lem-properties-spec-dpas}\hspace{-0.2cm}{\bf .}
\hspace{0.05cm}
Let $M \in \fcm$ be an $n$-dimensional with $k =[\frac {n-1}{2} ] $
odd and with a fixed spin structure.
Let $m$ be a natural number such that $m \equiv r  \mod (n) .$
Assume that $f_bv_\epsilon \in {\cal E}(\lambda , D_{Pas} ).$
Then

\vspace{0.1cm} \noindent
{\bf a)} $f_{-b}v_{-\epsilon } \in {\cal E}(\lambda , D_{Pas}) ,$

\vspace{0.1cm} \noindent
{\bf b)} $\dim {\cal E}(2\pi (m) , D_{Pas}) = A_{r}^{+}$ and 
$\dim {\cal E}(2\pi (m+\frac{1}{2}) , D_{Pas}) = A_{r}^{-}.$
\end{lemma}
\noindent
{\bf Proof. a)} If the spin structure on $M$ is
$\alpha _{+},$ then $b=le_n$ for some $l\in \Z ,$
and, from the equivariance of $f_bv_\epsilon,$ it follows that
$$ l \equiv \fracsmall {\mu _\epsilon }{2} +c(k)n \mod (n) .$$
Hence
$$ -l \equiv \fracsmall {\mu _{-\epsilon }  }{2} +c(k)n \mod (n) $$
and $f_{-b}v_{-\epsilon } $ is $g$-equivariant.
By the assumption that $k$ is odd,
$\nu (-\epsilon ) =-\nu (\epsilon ).$
According to Lemma ,
$f_{-b}v_{-\epsilon } \in {\cal E}(\lambda , D_{Pas}) .$

\vspace{0.2cm}
If the spin structure is $\alpha _{-},$ then
we use the congruence
$$ l \equiv \fracsmall {\mu _\epsilon }{2} +c(k)n +k \mod (n) .$$
Since
$\mu _{-\epsilon } =-\mu _\epsilon ,$
$c(k)n \equiv -c(k) \mod (n) ,$ and
$-k-1 \equiv k \mod (n) $   we have
$$ -l-1 \equiv \fracsmall {\mu _{-\epsilon } }{2} +c(k)n +k \mod (n) $$
and consequently $f_{-b}v_{-\epsilon }$ is
$g$-equivariant.

\vspace{0.2cm}
\noindent
{\bf b)}
If the spin structure is $\alpha _{+},$
then using Proposition \ref{prop-spectrum-Dpas} and a),
we get
$$
\dim {\cal E}(2\pi (m) , D_{Pas}) =
\dim {\cal E}(2\pi (r) , D_{Pas}) $$
$$ = 2 \card \{ \epsilon \in {\cal D}_{+} : \,
\fracsmall {\mu _\epsilon }{2} +c(k)n \equiv r \mod (n) \} =
A_{r}^{+} .$$
In the case of $\alpha _{-}$ 
we get 
$$\dim {\cal E}(2\pi (m + \frac{1}{2}) , D_{Pas}) = 
\dim {\cal E}(2\pi (r + \frac{1}{2}) , D_{Pas})$$
$$
= 2 \card \{ \epsilon \in {\cal D}_{+} :
 \frac {\mu _\epsilon }{2} +c(k)n + k
 \equiv  r \; \mod (n)\} = A_{r}^{-}.$$
\qed
\vskip 5mm
\noindent
\noindent {\bf Proof of Theorem \ref{thm-calcul-eta}.} 
\vspace{0.2cm}
 We
shall modify of
a proof of Lemma 5.5 from \cite{Pf}. 
 {\bf a)}
Let 
$m \in \Z $,
and let
${\cal S}_r =\{ 2\pi (m):\, m \equiv r$ mod $(n) \} .$
It is clear that ${\cal S}_r$ are disjoint and
${\cal S}_{Pas} \subseteq \Cup _{r=0}^{n-1}{\cal S}_r$.
Since $D_S$ has symmetric spectrum, 
$$ \eta _M(z) =
\sum _{\lambda \in {\cal S}_{Pas} }
\frac { \sgn (\lambda )}{\vert \lambda \vert ^z }
\dim {\cal E} (\lambda , D_{Pas } )$$
for ${\rm Re} (z)$ sufficiently large.
By Lemma \ref{lem-properties-spec-dpas} b),
$\dim {\cal E}(\lambda ,D_{Pas} ) =A_{r}^{+} $ for
$\lambda \in {\cal S}_{r}.$
If $A_{0}^{+}\neq 0$, then
the eigenvalue $\lambda \in {\cal S}_{0}$ occur together
with $-\lambda $ with the same multiplicity $A_{0}^{+}$
so that  
$\sum _{\lambda \in {\cal S}_0 -\{ 0\} }
\frac {A_{0}^{+} \, \sgn (\lambda ) } {\vert \lambda \vert ^z} =0 $
and, for ${\rm Re }(z)$ sufficiently big,
$$ \eta _M(z)=
\sum _{r =1 }^{n-1}
\sum _{m=-\infty } ^{\infty }
\frac {A_{r}^{+} \, \sgn (2\pi (mn+r) )}{ \vert 2\pi (mn+r)\vert ^z } =
\sum _{r =1 }^{n-1}
\sum _{m=-\infty } ^{\infty }
\frac {A_{r}^{+}\, \sgn \left (2\pi n  \left ( m+ \frac {r}{n} \right )
\right )  }
{\left \vert 2\pi n \left ( m+ \frac {r}{n} \right )  \right \vert ^z }$$
$$= \sum _{r = 1 }^{n-1}  \frac{A_{r}^{+}}{\mid 2\pi n \mid ^z} \left ( \sum_{m=0} ^{\infty}
\frac {1}{(m+\frac{r}{n})^z} - \sum_{m=0} ^{\infty} \frac{1}{(m+1-\frac{r}{n})^z} \right ).$$
The last two series are known as generalized zeta functions (cf. \cite{ww}).
They have meromorphic
extensions on $\C$ without poles in $z = 0.$ Let $\zeta(z,a)$ denote the function defined
by $\sum_{m=0} ^{\infty} \frac {1}{(m+\frac{r}{n})^z}$ for ${\rm Re }(z)$ sufficiently big. 
One gets for the extension:  $\zeta(0,a) = \frac{1}{2} - \frac{r}{n}.$ Hence  
$$ \eta _{M}(0) =
\sum _{r =1 }^{n-1}
A_{r}^{+} \left (1- \frac {2r}{n} \right ).$$
\vspace{0.2cm}
\noindent
{\bf b)}
We use similar arguments as those given in the proof of a).
Now the component ${\cal S}_0$ is not symmetric so
that we do not remove $r=0$ from the formula describing
$\eta _M(z).$
The equality
$$ 2\pi \left ( mn +r +\frac {1}{2} \right ) =
2\pi n \left ( m +\frac {2r+1}{2n} \right ) $$
and the above considerations
implies that
$$ \eta _M(0) =
\sum _{r =0 }^{n-1} A_{r}^{-} \left (1- \frac{2r+1}{n} \right ) . $$
This finishes the proof of Theorem
\ref{thm-calcul-eta}. \qed
\vskip 1mm
\noindent
We have.
\begin{corollary}
Let $n$ be a prime number greater than 3 such that $n+1$ is divisible by $4$
and let $M^{n}\in \fcm$ be a flat manifold with a fixed spin structure.
Then $\eta_{M^{n}} \in \Z.$
\end{corollary}
{\bf Proof:} Let $l = \frac{n+1}{4}.$ It is known (cf. \cite[chapter 9]{ston} that
$2^{s}$ copies of $M^{n}$ is a boundary of a spin manifold $W^{n+1}$ for some $s\in \N.$
By \cite[Theorem 4.2]{aps}
$$\int_{W^{n+1}}\hat{A}_{l}(p) - \frac{2^{s}\eta_{M^{n}}}{2} \in \Z,$$
where $\hat{A_{l}}$ is the $l$-th $\hat{A}$-polynomial on Pontriagin classes.
By \cite{hir} $\int_{W^{n+1}}\hat{A}_{l}(p)$ can be written as
$\frac{C_{W^{n+1}}}{q_1\dots q_r},$ where $C_{W^{n+1}}\in \Z$ and where
$q_1, \dots, q_r \in \{2,3,... , n-1\}$ are prime numbers.
From Theorem \ref{thm-calcul-eta} $\eta_{M^{n}} = \frac{C_{M^{n}}}{n},$ for some $C_{M^{n}} \in \Z.$
Since $\frac{C_{W^{n+1}}}{q_1\dots q_r} - 2^{s-1}\eta_{M^{n}} \in \Z$
we have $\frac{2^{s-1}q_1\dots q_r C_{M^{n}}}{n} \in \Z.$ Hence $\eta_{M^{n}} \in \Z.$
\qed
\begin{corollary} Let $n$ be an odd number.
If $M^{n}\in \fcm,$ then $d = \eta_{M^{n},\alpha^{+}} - \eta_{M^{n},\alpha^{-}} \in 2\Z.$
\end{corollary}
{\bf Proof:} By \cite[Theorem 1.1]{dahl}, $d\in \frac{1}{2}\Z.$ From the definition (cf. page 2)  
all $A_{r}^{\pm}$ belongs to $2\Z.$ 
Hence $d = \frac{2C}{n}$ for some $C\in \Z.$ Summing up $\frac{dn}{2} \in \Z$
and $d \in 2\Z.$ 
\qed
\begin{examp}\label{ex-eta-in-dim-7}\hspace{-0.2cm}{\bf .} \hspace{0.05cm}
\begin{rm}
We calculate $\eta_{M^{7},\alpha _{+}},$ where$M^{7}\in \fcm .$
Since $\frac{k(k+1)}{2} = 6$ is even our equation is
$$ r \equiv \frac {\mu _\epsilon }{2}  \mod (7) .$$
The values of $\frac {\mu _\epsilon }{2}$ and $r$ for
$\epsilon \in {\cal D}_{+} $ are given in the following table.

\vspace{0.2cm} \noindent
\begin{center}
\begin{tabular}{|l|r|l|}
\hline
$\epsilon $ \hspace{0.3cm} & $\frac {\mu _\epsilon }{2}  $  & $r  $ \\
\hline
$(1,1,1)$   & $3$    & $3$ \\
$(1,-1,-1)$   & $-2  $ & $5 $\\
$(-1,1,-1)$   & $-1 $ &  $ 6$\\
$(-1,-1,1)$   & $0  $  & $0$ \\
\hline
\end{tabular}
\end{center}
\vspace{0.2cm} \noindent
It follows that  $A_{j}^{+}=2,$ for $j=0,3,5,6 $
and $A_{j}^{+}=0$ for other values of $j.$
By Theorem \ref{thm-calcul-eta},
$$ \eta _M^{7}(0)=
2 \left [
\left (1-\frac {6}{7} \right )+
\left (1-\frac {10}{7} \right )+
\left (1-\frac {12}{7} \right  ) \right ] =
-2 .$$
\end{rm}
\end{examp}

\section{Harmonic spinors on $\fcm $-manifolds}\label{sec-harmonin-fcm}

{\it A harmonic spinor} on a closed spin manifold $M$
is an element of the kernel of the Dirac operator on $M.$

\vspace{0.2cm}
\noindent
{\bf Proof of Proposition 1.  a)}
We have
$$ g v_\epsilon =(-1)^{\frac {k(k+1)}{2} } \alpha v_\epsilon =
(-1)^{\frac {k(k+1)}{2} } e^{\frac {2\pi i}{2n}\mu _\epsilon } v_\epsilon
.$$
First consider the case when $k(k+1)/2$ is even.
Then $g v_\epsilon = v_\epsilon $ if and only if
$$ \mu _\epsilon \equiv 0 \mod (2n) $$
and $k=4k_0 +3$ or $k=4k_0 .$
Let $\delta _4 $ denote the sequence
$1,-1,-1,1 .$
If $k=4k_0 +3 $ and
$$ \epsilon = (-1,-1,1, \underbrace {\delta _4,..., \delta _4}_
{k_0 \; {\rm times}} ) , $$
then $\epsilon $
belongs to $\{-1,1 \} ^k $ and $\mu _\epsilon = 0. $
If $k=4k_0 $ and
$$ \epsilon = (\underbrace {\delta _4,..., \delta _4}_
{k_0 \; {\rm times}} ) ,$$
then $\epsilon $
belongs to $\{-1,1 \} ^k $ and $\mu _\epsilon = 0. $
In particular, ${\mathfrak h}(M)>0. $

\vspace{0.2cm}
Now assume that $k>2 $ and $k(k+1)/2$ is odd.
Then $g v_\epsilon = v_\epsilon $ if and only if
$$ \mu _\epsilon \equiv n \mod (2n) $$
and $k=4k_0 +1$ or $k=4k_0 +2 .$
If $k=4k_0 +1 $ consider
$$ \epsilon = (1,-1,1, \underbrace {\delta _4,..., \delta _4}_
{k_0 -1 \; {\rm times}},1,1 ) $$
and if $k=4k_0 +2$ consider
$$ \epsilon = (-1,1-1,1,\underbrace {\delta _4,..., \delta _4}_
{k_0 -1 \; {\rm times}}  ,1,1) .$$
In both cases  $\epsilon \in \{-1,1 \} ^k $ and $\mu _\epsilon = n. $
It is easily seen that the equation
$\mu _\epsilon \equiv 0 \mod (2n) $ have no solutions for
$k=1$ or $2$.

\vspace{0.2cm}
\noindent {\bf b)}
Since $g^n=-id$, the equation
$g v = v $ have only one solution $v=0.$  \qed

\vspace{0.2cm}
It is easy to see that the equality
$\alpha _{+} v_\epsilon = v_\epsilon $ implies
$\alpha _{+} v_{-\epsilon }= v_ {-\epsilon }.  $
Using this and the arguments given in the proof of
Proposition 1 we have.

\begin{corollary}
\hspace{0.05cm}
If $M\in \fcm $ and $\dim M= 2k+1,$ then
$${\mathfrak h}(M, \alpha _{+})=
2 \card \{ \epsilon \in {\cal D}_{+} :\,
\frac {\mu _\epsilon }{2} +c(k)n \equiv 0 \mod (2k+1) \}.$$
\end{corollary}


\noindent 
Institute of Mathematics\\
University of Gda\'nsk\\
ul.Wita Stwosza 57\\
80 - 952 Gda\'nsk\\
Poland\\
E - mail : msa @ delta.math.univ.gda.pl\\
E - mail : matas @ paula.univ.gda.pl


\end{document}